\documentclass[11pt, a4paper]{article}

\usepackage[utf8]{inputenc}
\usepackage[english]{babel}
\usepackage[T1]{fontenc}

\usepackage{amsmath}
\usepackage{amssymb}
\usepackage{amsthm}
\usepackage{mathabx}
\usepackage{tikz}
\usepackage{tikz-cd}
\usepackage{todonotes}
\usepackage{hyperref}
\usepackage{cleveref}
\usepackage{footmisc}
\usepackage{bbm}
\usepackage{afterpage}
\usepackage{mathrsfs}
\usepackage{sseq}
\usepackage{extarrows}
\usepackage{enumitem}
\usepackage{stmaryrd}
\usepackage{pdflscape}

\newcommand{\alg}{\mathrm{alg}}

\newcommand{\chrom}{{\mathrm{chrom}}}

\DeclareMathOperator{\Spec}{Spec}
\DeclareMathOperator{\Spf}{Spf}

\newcommand{\sm}{{\mathrm{sm}}}

\newcommand{\Ell}{\mathrm{Ell}}
\newcommand{\et}{\mathrm{\acute{e}t}}

\newcommand{\ord}{{\mathrm{ord}}}
\newcommand{\supers}{{\mathrm{ss}}}

\DeclareMathOperator{\KU}{KU}

\DeclareMathOperator{\Tmf}{Tmf}
\DeclareMathOperator{\tmf}{tmf}

\newcommand{\Alg}{{\mathrm{Alg}}}

\newcommand{\CAlg}{{\mathrm{CAlg}}}

\DeclareMathOperator{\Fun}{Fun}

\newcommand{\Sp}{\mathrm{Sp}}
\newcommand{\Spc}{\mathcal{S}}
\renewcommand{\top}{\mathrm{top}}

\renewcommand{\lim}{{\mathrm{lim}\,}}

\newcommand{\Hom}{\mathrm{Hom}}

\DeclareMathOperator{\Map}{Map}

\DeclareMathOperator{\Tw}{Tw}

\newcommand{\op}{\mathrm{op}}

\newcommand{\Z}{\mathbf{Z}}

\newcommand{\A}{\mathcal{A}}

\newcommand{\B}{\mathcal{B}}

\newcommand{\C}{\mathcal{C}}

\newcommand{\D}{\mathcal{D}}

\newcommand{\E}{\mathbf{E}}
\newcommand{\EE}{\mathcal{E}}

\newcommand{\F}{\mathbf{F}}
\newcommand{\FF}{\mathcal{F}}

\newcommand{\G}{\mathbf{G}}

\newcommand{\h}{\mathrm{h}}

\newcommand{\K}{\mathrm{K}}

\newcommand{\M}{\mathcal{M}}

\renewcommand{\O}{\mathscr{O}}

\newcommand{\Q}{\mathbf{Q}}
\newcommand{\QQ}{\mathcal{Q}}

\newcommand{\U}{\mathcal{U}}

\renewcommand{\Z}{\mathcal{Z}}

\newcommand{\al}{\alpha}
\newcommand{\be}{\beta}
\newcommand{\ga}{\gamma}
\newcommand{\Ga}{\Gamma}

\newenvironment{proofclaim}[1]
{\emph{#1.}}
{\hfill $\blacksquare$\break}

\usepackage{mathtools}

\theoremstyle{theorem}%\numberwithin{equation}{section}
\newtheorem{theorem}[equation]{Theorem}
\newtheorem{theoremalph}{Theorem}

\newtheorem{prop}[equation]{Proposition}
\newtheorem{lemma}[equation]{Lemma}

\theoremstyle{definition}%\numberwithin{equation}{section}
\newtheorem{mydef}[equation]{Definition}

\theoremstyle{remark}%\numberwithin{equation}{section}

\newtheorem{claim}[equation]{Claim}

\newtheorem{remark}[equation]{Remark}

\begin{document}
\title{Elliptic cohomology is unique up to homotopy}
\author{J.M. Davies\footnote{\url{j.m.davies@uu.nl}}}
\maketitle

\begin{abstract}
Homotopy theory folklore tells us that the sheaf defining the cohomology theory $\Tmf$ of topological modular forms is unique up to homotopy. Here we provide a proof of this fact, although we claim no originality for the statement. This retroactively reconciles all previous constructions of $\Tmf$.
\end{abstract}

\setcounter{tocdepth}{1}
\tableofcontents

\addcontentsline{toc}{section}{Introduction}
\section*{Introduction}
At present, there is a singular definition of the extraordinary cohomology theory $\Tmf$ of \emph{topological modular forms}, and that is as the global sections $\Ga(\M_\Ell,\O^\top)$. Here, and elsewhere in this note, $\M_\Ell$ is the compactification of the moduli stack of elliptic curves and $\O^\top$ is the Goerss--Hopkins--Miller--Lurie sheaf of $\E_\infty$-rings; we recommend \cite{handbooktmf} for more background. The first published construction of this sheaf is due to Behrens \cite[\textsection12]{tmfbook}, and is based on the unpublished work of Goerss, Hopkins, and Miller.\\

These days one can find variations of this construction in the literature, some using logarithmic geometry on the moduli stack $\M_\Ell$ \cite{hilllawson}, and others using oriented derived elliptic curves \cite{ec2name}. This leads us to a fundamental question:
\begin{center}
\emph{Are these different constructions in any way compatible?}
\end{center}
This question would be easily answered if one could show that $\O^\top$ is the ``unique'' sheaf with a certain property, assuming this property holds for each competing construction. Moreover, such uniqueness could also be used to further justify the significance of $\O^\top$; we leave it to the reader to come to her own conclusions on that front. The property we would like to consider is that of an \emph{elliptic cohomology theory}.

\begin{mydef}\label{genectheory}
Let $E$ be a generalised elliptic curve over a ring $R$ with irreducible geometric fibres, which is equivalent data to a morphism of stacks $\Spec R\to \M_\Ell$. We say that a homotopy commutative ring spectrum $\EE$ is an \emph{elliptic cohomology theory}\footnote{These were previously called \emph{generalised elliptic cohomology theories} in \cite[\textsection2.2]{adamsontmf}, and other variations on this theme can be found elsewhere; see \cite[\textsection6]{marktmf} or \cite[Df.1.2]{lurieecsurveyname}, for example.} for $E$ (or $\Spec R\to \M_\Ell$) if we have the following data:
\begin{enumerate}
\item $\EE$ is weakly 2-periodic, meaning the homotopy group $\pi_2 \EE$ is a projective $\pi_0 \EE$-module of rank one and for every integer $n$ the canonical map of $\pi_0\EE$-modules
\[\pi_2 \EE\underset{\pi_0\EE}{\otimes} \pi_n \EE\to \pi_{n+2}\EE\]
is an isomorphism; see \cite[\textsection4.1]{ec2name}; 
\item The groups $\pi_k \EE$ vanish for all odd integers $k$, so in particular $\EE$ is complex orientable; 
\item There is a chosen isomorphism of rings $\pi_0 \EE\simeq R$; and
\item There is a chosen isomorphism of formal groups $\widehat{E}\simeq \widehat{\G}_\EE$ over $R$, between the formal group of $E$ and the classical Quillen formal group of $\EE$; see \cite[\textsection4]{ec2name}.
\end{enumerate}
We say a collection of such $\EE$ is \emph{natural} if the isomorphisms in parts 3-4 above are natural with respect to some subcategory of affine schemes over $\M_\Ell$.
\end{mydef}

The following is a simple uniqueness statement for $\O^\top$ as a functor valued in homotopy commutative ring spectra.

\begin{prop}\label{weakunique}
The functor $\h\O^\top\colon \U^\op\to \CAlg(\h\Sp)$, from the small affine \'{e}tale site $\U$ of $\M_\Ell$ to the 1-category of homotopy commutative ring spectra, is uniquely defined up to isomorphism by the property that it defines natural elliptic cohomology theories on $\U$.
\end{prop}

The proof of the above statement follows from the fact that each section $\O^\top(R)$ is Landweber exact; see \cite[Rmk.1.6]{marktmf} for the statement and \cite[\textsection A]{adamsontmf} for more discussion and a proof. A remarkable fact about $\O^\top$ is that the property that it defines a natural elliptic cohomology theory actually characterises this sheaf with values in the $\infty$-category $\CAlg$ of $\E_\infty$-rings. The following is stated (without proof) in \cite[Th.1.1]{lurieecsurveyname} and \cite[Th.1.2]{bourbakigoerss}.

\begin{theoremalph}\label{tha}
The sheaf of $\E_\infty$-rings $\O^\top$ on the small \'{e}tale site of $\M_\Ell$ is uniquely defined up to homotopy by the property that it defines natural elliptic cohomology theories on the small affine \'{e}tale site of $\M_\Ell$. The same holds for the restriction $\O^\top_\sm$ of $\O^\top$ to the small \'{e}tale site of $\M_\Ell^\sm$, the moduli stack of \emph{smooth elliptic curves}.
\end{theoremalph}

The difference between the \Cref{weakunique} and \Cref{tha} is two-fold, as \Cref{tha} says not only that $\O^\top$ is uniquely defined up to homotopy as a sheaf of $\E_\infty$-rings (as opposed to a diagram in $\CAlg(\h\Sp)$), but that statement is made in an $\infty$-category (as opposed to the 1-category $\CAlg(\h\Sp)$).\\

The utility of \Cref{tha} is evident. For example, it retroactively shows that the various constructions of $\O^\top$ found in \cite{marktmf}, \cite{hilllawson}, and \cite[\textsection2]{adamsontmf} (and also \cite[\textsection7]{ec2name} and \cite[\textsection2.3]{luriestheorem} over the moduli stack of smooth elliptic curves) all agree up to homotopy. Importantly, \Cref{tha} constructs noncanonical (see \Cref{higherstuff}) equivalences of $\E_\infty$-rings between all available definitions of $\Tmf$; a conclusion which does \textbf{not} follow directly from \Cref{weakunique}. The author also finds the proof long and complex enough to warrant a publicly available write-up.

\subsection*{Outline}
To prove \Cref{tha}, we will first reduce the question to one of the connectedness of a certain moduli space; see \textsection\ref{statementsection}. In \textsection\ref{naturalitysection}, we formulate and prove a statement about spaces of natural transformations which we will often use; we suggest that the reader initially skips this section, and only returns when she deems it necessary. Our proof of \Cref{tha} (reformulated as \Cref{mainguy}) occurs in \textsection\ref{proofsection}, and follows Behrens' construction of $\Tmf$ rather closely: first we work with the separate chromatic layers, before gluing things together in both a transchromatic sense and then an arithmetic sense. The $K(1)$-local case in this section requires a statement about $p$-adic Adams operations on $p$-adic $K$-theory, which is the focus of \textsection\ref{finalcompsection}.

\subsection*{Conventions}
We will assume that the reader is familiar with the language of $\infty$-categories as well as the techniques used in the construction of $\O^\top$ as described by Behrens \cite{marktmf}. Furthermore, we advise the reader keeps a copy of \emph{idem} in her vicinity. Write $\Map_\C(X,Y)$ for the mapping space between any two objects $X,Y$ in an $\infty$-category $\C$, and $\Hom_{\h\C}(X,Y)$ for its zeroth truncation, or equivalently, for $\pi_0\Map_\C(X,Y)$. Let us also suppress the notation distinguishing a 1-category from its nerve, considered as an $\infty$-category, and the same for a 2-category, such as the small \'{e}tale site over $\M_\Ell$; see \cite[\href{https://kerodon.net/tag/007J}{Tag 007J}]{kerodon}.

\subsection*{Acknowledgement}
I owe a huge thank you to Lennart Meier, for teaching me everything I know about $\O^\top$.

%%%%%%%%%%%%%%%%%%%%%%%%%%%%%%%%%%%%%%%%%%%%%%%%%%%%

\section{A reformulation}\label{statementsection}
Let us now make a statement to help us prove \Cref{tha}.

\begin{theorem}\label{mainguy}
Write $\U$ (resp. $\U_\sm$) for the (2-) category of affine schemes with \'{e}tale maps to $\M_\Ell$ (resp. $\M_\Ell^\sm$). Then the spaces
\[\Z=\Fun\left(\U^\op, \CAlg\right)\underset{\Fun\left(\U^\op, \CAlg(\h\Sp)\right)}{\times} \{\h\O^\top\}\]
\[\Z^\sm=\Fun\left(\U^\op, \CAlg\right)\underset{\Fun\left(\U_\sm^\op, \CAlg(\h\Sp)\right)}{\times} \{\h\O^\top\}\]
are connected.
\end{theorem}

\begin{remark}\label{higherstuff}
As mentioned in \cite[Rmk.7.0.2]{ec2name}, the moduli space $\Z^\sm$ is \textbf{not} contractible. In other words, \Cref{mainguy} states that $\O^\top$ is unique as a $\CAlg$-valued presheaf of elliptic cohomology theories on $\U^\sm$ only up to homotopy, and \textbf{not} up to contractible choice. We would like to guide the reader to an explanation for this fact given by Tyler Lawson on \url{mathoverflow.net}; see \cite{tylersanswer}.
\end{remark}

\begin{proof}[Proof of \Cref{tha} from \Cref{mainguy}]
The $\infty$-category of sheaves of $\E_\infty$-rings on the \'{e}tale site of $\M_\Ell$ is equivalent, by restriction and right Kan extension, to the $\infty$-category of sheaves of $\E_\infty$-rings on the affine \'{e}tale site of $\M_\Ell$; see \cite[Lm.2.6]{adamsontmf} for a similar argument following the ``comparison lemma'' of \cite[Lm.C.3]{hoyoisappendixtimesaver}. Note that the latter is an $\infty$-subcategory of $\Fun(\U^\op,\CAlg)$, and that if a functor $F\colon \U^\op\to\CAlg$ defines an natural elliptic cohomology theory and there is an equivalence $F\simeq G$, then $G$ also defines a natural elliptic cohomology theory. These two observations show that it suffices to prove the space $\Z'$ is connected, where $\Z'$ is the component of $\Fun\left(\U^\op, \CAlg\right)^\simeq$ spanned by those functors which define natural elliptic cohomology theories. There is a map $\Z\to \Z'$ as both $\O^\top$ and any presheaf of $\E_\infty$-rings equivalent to $\O^\top$ as a diagram of homotopy commutative ring spectra, defines a natural elliptic cohomology theory. The map $\Z\to \Z'$ induces an equivalence on $\pi_0$ as \cite[Pr.A.1]{adamsontmf} (also see \cite[Rmk.1.6]{marktmf}) states that any functor $\U^\op\to \CAlg(\h\Sp)$ which defines an elliptic cohomology theory is isomorphic to $\h\O^\top$. \Cref{mainguy} then implies that the moduli space $\Z'$ is also connected. The same argument can be made for $\Z^\sm$.
\end{proof}

\begin{remark}
Write $\U_\Q$ for the small affine \'{e}tale site of $\M_\Ell\times\Spec \Q$ and for each prime $p$ write $\U_p$ for the small affine \'{e}tale site of $\M_\Ell\times\Spf \mathbf{Z}_p$. The construction of $\O^\top$ as found in \cite{marktmf} for example proceeds first with a rational construction $\O^\top_\Q$ over $\U_\Q$, and a $p$-complete construction $\O^\top_p$ over $\U_p$. The methods of our proof for \Cref{mainguy} show that the moduli spaces $\Z_\Q$ and $\Z_p$, of realisations of $\h\O^\top_\Q$ and $\h\O^\top_p$ over these aforementioned sites, are also connected. This means that analogs of \Cref{tha} also holds for both $\O^\top_\Q$ and $\O^\top_p$. The same holds for the $p$-completed and rational version of $\O^\top_\sm$ for similar reasons. Moreover, following the ``arithmetic compatibility'' discussed in the proof of \Cref{mainguy}, it follows that the localisations $\O^\top[\mathcal{P}^{-1}]$ and $\O_\sm^\top[\mathcal{P}^{-1}]$ satisfy their own version of \Cref{tha}, where $\mathcal{P}$ is any set of primes.
\end{remark}

The following is a short remark on the homotopy groups of elliptic cohomology theories which will be important later.

\begin{remark}\label{higherhomotopy}
Let $\EE$ be an $\E_\infty$-elliptic cohomology theory for some $E\colon\Spec R\to \M_\Ell$. It follows that there is a natural isomorphism $\pi_{2k}\EE\simeq \omega_E^{\otimes k}$ for all integers $k$, where $\omega_E$ is the dualising line for the formal group $\widehat{E}$; see \cite[\textsection4.2]{ec2name}. Indeed, as the odd homotopy groups of $\EE$ vanish, we see $\EE$ possesses a complex orientation which yields the classical Quillen formal group $\widehat{\G}_\EE^{\QQ_0}$ over $\pi_0 \EE$; see \cite[Ex.4.1.2]{ec2name}, for example. From this we see $\EE$ is \emph{complex periodic}, meaning it has a complex orientation and is weakly 2-periodic (see \cite[\textsection4.1]{ec2name}), and \cite[Ex.4.2.19]{ec2name} then implies that $\pi_2\EE$ is isomorphic to the dualising line for the formal group $\widehat{\G}_\EE^{\QQ_0}$. Part 4 of \Cref{genectheory} states that $\pi_2\EE$ is naturally isomorphic to $\omega_E$, and part 1 gives us the claim above.
\end{remark}

%%%%%%%%%%%%%%%%%%%%%%%%%%%%%%%%%%%%%%%%%%%%%%%%%%%%

\section{Spaces of natural transformations}\label{naturalitysection}
To prove \Cref{mainguy}, we will show that any two functors $\O$ and $\O'$ in $\Z$ can be connected by a path in $\Z$. In particular, we would like effective tools for studying spaces of natural transformations between functors of $\infty$-categories. The following is known to experts, and a model categorical interpretation can be found in \cite{DKS89}.

\begin{prop}\label{naturality}
Let $\C,\D$ be $\infty$-categories and $F,G\colon \C\to \D$ be functors. Suppose that for all objects $X,Y$ in $\C$ the mapping space $\Map_\D(FX,GY)$ is discrete, meaning the natural map
\[\Map_\D(FX,GX)\to \Hom_{\h\D}(\h FX, \h GY)\]
is an equivalence of spaces. Then the mapping space $\Map_{\Fun(\C,\D)}(F,G)$ is also discrete, so the natural map
\[\Map_{\Fun(\C,\D)}(F,G)\to \Hom_{\h\Fun(\C,\D)}(F,G)\simeq \Hom_{\Fun(\C,\h\D)}(\h F,\h G)\]
is an equivalence of spaces, where an $\h$ before a functor denotes post composition with the natural map $\D\to\h\D$; the unit of the homotopy category-nerve adjunction of \cite[Pr.1.2.3.1]{httname}.
\end{prop}

\begin{proof}
By \cite[Pr.5.1]{davidrunethomas}, the space of natural transformations from $F$ to $G$ is naturally equivalent to the limit of the diagram
\begin{equation}\label{endone}\Tw(\C)^\op\xrightarrow{H}\C^\op\times \C\xrightarrow{F^\op\times G} \D^\op\times \D\xrightarrow{\Map_\D(-,-)} \Spc,\end{equation}
where $\Tw(\C)$ is the \emph{twisted arrow category} of $\C$ (see \cite[Df.2.2]{davidrunethomas}), and the $\Tw(\C)\to\C\times\C^\op$ is the natural right fibration (see \emph{idem}).\footnote{We will stick to the notation and conventions of \cite{davidrunethomas}, which is a particular choice out of a possible two; see \cite[Wrn.2.4]{davidrunethomas}.} The limit of (\ref{endone}) is by definition the \emph{end} of the composition $\C^\op\times\C\to \Spc$. Consider the following not a priori commutative diagram of $\infty$-categories:
\begin{equation}\label{littlethreesquarediagram}\begin{tikzcd}
{\Tw(\C)^\op}\ar[d]\ar[r, "{T}"]	&	{\C^\op\times\C}\ar[d]\ar[rrrrrd, dashed]\ar[rr, "{F^\op\times G}"]	&&	{\D^\op\times\D}\ar[d]\ar[rrr, "{\Map_\D(-,-)=M}"]	&&&	{\Spc}	\\
{\h\Tw(\C)^\op}\ar[r, "{T'}", swap]	&	{\h\C^\op\times\h\C}\ar[rr, "{\h F^\op\times\h G}", swap]	&&	{\h\D^\op\times\h\D}\ar[rrr, "{\Hom_{\h\D}(-,-)=H}", swap]	&&&	{\Spc_{\leq0}}\ar[u]	\\
\end{tikzcd}\end{equation}
Above, the vertical functors are the obvious ones, hence the left and middle squares commute. Our hypotheses dictate that the dashed arrow above exists, which we will now denote by $P$, such that the top-right and bottom left triangles commutes. As the inclusion of $\infty$-subcategories $\Spc_{\leq 0}\subseteq \Spc$, from the $\infty$-category of discrete spaces, preserves limits, we note it suffices to compute the limit of (\ref{endone}) as the limit of $P\circ T$ inside $\Spc_{\leq0}$. As this limit lands in $\Spc_{\leq 0}$, which is equivalent to the nerve of the 1-category of sets, we see the limit of $P\circ T$ can be calculated as the limit of the lower-horizontal composition of (\ref{littlethreesquarediagram}). We then obtain the following natural equivalences, twice employing \cite[Pr.5.1]{davidrunethomas}, first for general $\infty$-categories, and again in the classical 1-categorical case:
\[\Map_{\Fun(\C,\D)}(F,G)\simeq \underset{\Tw(\C)}{\lim} M (F^\op\times G) T\simeq \underset{\Tw(\C)}{\lim} PT\]
\[\simeq \underset{\Tw(\h\C)}{\lim} H (\h F^\op\times \h G)T'\simeq \Hom_{\Fun(\h\C,\h\D)}(\h F,\h G).\]
The final (discrete) space above is naturally equivalent to $\Hom_{\Fun(\C,\h\D)}(\h F, \h G)$ from the natural equivalence of $\infty$-categories $\Fun(\C,\h\D)\simeq\Fun(\h\C,\h\D)$.
\end{proof}

%%%%%%%%%%%%%%%%%%%%%%%%%%%%%%%%%%%%%%%%%%%%%%%%%%%%

\section{The proof of \Cref{mainguy}}\label{proofsection}
Let $\O\colon \U^\op\to \CAlg$ be an object of $\Z$, hence it comes equipped with an equivalence $\h\phi\colon\h\O^\top\to \h\O$ of functors $\U^\op\to \CAlg(\h\Sp)$. To see $\Z$ is connected, it suffices to show $\h\phi$ can be lifted to an equivalence $\phi\colon \O^\top\to \O$ of presheaves of $\E_\infty$-rings on $\U$. Fix such an $\h\phi$ for the remainder of this proof. Let us work section-wise, so we also fix an object $\Spec R\to \M_\Ell$ inside $\U$, and write
\[\h\phi\colon \EE^\top=\O^\top(R)\to \O(R)=\EE\]
for the given natural equivalence of homotopy commutative ring spectra. To naturally lift this map to one of $\E_\infty$-rings we will work through the layers of chromatic homotopy theory. This means we will first work $K(2)$-locally, $K(1)$-locally, and then $K(0)$-locally, where $K(n)$ denotes the \emph{$n$th Morava $K$-theory spectrum at a prime $p$}, before gluing these cases together with a $p$-complete statement followed by an arithmetic statement.

%%%%%%%%%%%%%

\paragraph{($K(2)$-local case)}
Writing $\widehat{(-)}$ for base-change over $\Spf\mathbf{Z}_p$, we define $\Spf R^\supers\to \M_\Ell^\supers$ as the base-change of $\Spf \widehat{R}\to \widehat{\M}_\Ell$ over $\M_\Ell^\supers$, where the latter is the completion of $\widehat{\M}_\Ell$ at the moduli stack $\M_{\Ell,\F_p}^\supers$ of supersingular elliptic curves over $\F_p$. This pullback $\Spf R^\supers$ is affine by \cite[Rmk.8.7]{marktmf}. Write $E^\supers$ for the elliptic curve defined by $\Spf R^\supers\to \M_\Ell^\supers$. Serre--Tate and Lubin--Tate theory yield another description of $R^\supers$. Indeed, as $\M^\supers_{\Ell,\F_p}$ is zero-dimensional and smooth over $\Spec \F_p$, it follows that $\Spec R^\supers/I$ is \'{e}tale over $\F_p$, where $I$ is the finitely generated ideal generating the topology on $R^\supers$. This implies $R^\supers/I$ splits as a finite product $\prod_i \kappa_i$ where each $\kappa_i$ is a finite field of characteristic $p$. This provides a splitting of $E_0$, the reduction of $E^\supers$ over $R/I$, into $E_0\simeq \coprod E_0^i$. Writing $R_i\simeq W(\kappa_i)\llbracket u_1\rrbracket$ for the universal deformation ring of the pair $(\kappa_i, \widehat{E}^i_0)$ with associated universal formal group $\widehat{E}^\supers_{R_i}$, we obtain a natural equivalence $R^\supers\simeq \prod_i R_i$ as $E^\supers\colon \Spf R^\supers\to \M_\Ell^\supers$ was \'{e}tale; see \cite[Cor.4.3]{marktmf}.\\

By \cite[Pr.4.4]{marktmf}, the $K(2)$-localisations $\EE^\top_{K(2)}$ and $\EE_{K(2)}$ are elliptic cohomology theories for $R^\supers$, and also split into products $\EE^\top_i$ and $\EE_i$. It follows from \cite[\textsection7]{gh04} (also see \cite[Rmk.5.0.5]{ec2name} or \cite[\textsection7]{pp}), that these $K(2)$-local $\E_\infty$-rings $\EE_{K(2)}^\top$ and $\EE_{K(2)}$ are naturally equivalent to the product of Lubin--Tate $\E_\infty$-rings associated to the formal groups $\widehat{E}^\supers_{\kappa_i}$ over the (finite and hence also) perfect fields $\kappa_i$. By \emph{idem}, we see that morphisms between these Lubin--Tate $\E_\infty$-rings are defined by the associated morphisms on the pairs $(\kappa_i, \widehat{E}^\supers_{\kappa_i})$. As $\h\phi_{K(2)}$ yields an equivalence on $\pi_0$ as well as an equivalence on associated Quillen formal groups, we see $\h\phi_{K(2)}$ lifts to a morphism $\phi_{K(2)}\colon \EE^\top_{K(2)}\to \EE_{K(2)}$ of $K(2)$-local $\E_\infty$-rings, which is unique up to contractible choice. This uniqueness allows us to use \Cref{naturality} to conclude that the collection of morphisms of $\E_\infty$-rings define a natural morphism $\phi_{K(2)}\colon \O^\top_{K(2)}\to \O_{K(2)}$ of presheaves of $\E_\infty$-rings on $\U$; here $(-)_{K(2)}$ denotes $K(2)$-localisation.

%%%%%%%%%%%%%

\paragraph{($K(1)$-local case)}
Consider the $K(1)$-localisation $\h\phi_{K(1)}\colon \EE^\top_{K(1)}=\O^\top_{K(1)}(R)\to \O_{K(1)}(R)=\EE_{K(1)}$ of the map $\h\phi$ of homotopy commutative ring spectra. Recall from \cite[\textsection6]{marktmf}, that the $p$-adic $K$-theory\footnote{Recall that for a spectrum $X$, one defines its \emph{$p$-adic $K$-theory} as the homotopy groups of the localisation $\K_\ast^\wedge X=\pi_\ast L_{K(1)}(X\otimes \KU_p)$, or equivalently $\pi_\ast ((X\otimes \KU_p)_p^\wedge)$.} of an $\E_\infty$-ring has the structure of a $\theta$-$\pi_\ast\KU_p$-algebra, functorially in maps of $\E_\infty$-rings. Let us write $\M_\Ell^\ord$ for the moduli of generalised elliptic curves over $p$-complete rings with ordinary reduction modulo $p$ (see \cite[p.3]{marktmf}, and $\M_\Ell^\ord(p^\infty)$ for the moduli stack of generalised elliptic curves $E$ over $p$-complete rings and level structures given by an isomorphism $\widehat{\G}_m\simeq \widehat{E}$ of formal groups.

\begin{claim}\label{claimvanishing}
The following facts holds for the $p$-adic $K$-theory of $\EE_{K(1)}^\top$ and $\EE_{K(1)}$:
\begin{enumerate}
\item Both are isomorphic in degree zero to the pullback of $\Spec R\to \M_\Ell$ with the composite
\[\M_\Ell^\ord(p^\infty)\to \M_\Ell^\ord\to \widehat{\M}_\Ell\to \M_\Ell;\]
\item Both are concentrated in even degrees;
\item Both are ind-\'{e}tale over $R^\ord$, the base-change of $\Spec R$ over $\M_\Ell^\ord$, in degree zero; and
\item Both degree zero components have vanishing higher continuous group cohomologies, with respect to its $\mathbf{Z}^\times_p$-action inherited from part 1.
\end{enumerate}
\end{claim}

\begin{proofclaim}{Proof of \Cref{claimvanishing}}
Part 1 is obtained from \cite[Pr.6.1]{marktmf} by base-change. Part 2 also follows from a graded version of \cite[Pr.6.1]{marktmf}. Parts 3 and 4 follow from part 1, as the map $\M_\Ell^\ord(p^\infty)\to \M_\Ell^\ord$ is not only ind-\'{e}tale, but a $\mathbf{Z}_p^\times$-torsor; see \cite[Lm.5.1]{marktmf}.
\end{proofclaim}

For another object $\Spec R'\to \M_\Ell$ inside $\U$, consider the map induced by the $p$-adic $K$-theory functor
\begin{equation}\label{firstthetaalgebradiscreteness}
\Map_{\CAlg_{K(1)}}(\EE^\top_{K(1)}, \EE_{K(1)}')\to \Hom_{\theta\Alg_{\K^\wedge_\ast}}(\K^\wedge_\ast \EE^\top_{K(1)}, \K^\wedge_\ast \EE_{K(1)}'),
\end{equation}
where $\EE_{K(1)}'=\O_{K(1)}(R')$. Combining \Cref{claimvanishing} and the fact that $R^\ord$ is smooth over $\mathbf{Z}_p$ with \cite[Lm.7.5]{marktmf} implies that all the Andr\'{e}--Quillen cohomology groups in \cite[Th.7.1]{marktmf} vanish. By \emph{ibid}, it then follows that the above map is an equivalence of spaces. Despite the fact that each $\h\phi_{K(1)}(R)$ is currently just a morphism of homotopy commutative ring spectra, \Cref{thetaalgebra} states that its zeroth $p$-adic $K$-theory is a morphism of $\theta$-algebras. As $\mathbf{Z}$-graded $p$-adic $K$-theory obtains a $\theta$-algebra structure from that in degree zero, the $p$-adic $K$-theory of $\h\phi_{K(1)}$ defines an element inside the codomain of (\ref{firstthetaalgebradiscreteness}) when $R'=R$. By \Cref{naturality}, we can therefore lift $\h\phi_{K(1)}\colon \h\O^\top_{K(1)}\to \h\O_{K(1)}$ to a morphism $\phi_{K(1)}\colon \O^\top_{K(1)}\to \O_{K(1)}$ of presheaves of $\E_\infty$-rings on $\U$. 

%%%%%%%%%%%%%

\paragraph{($K(0)$-local case)}
The Morava $K$-theory spectrum $K(0)$ is equivalent to $\Q$, the Eilenberg--Mac~Lane spectrum of the rational numbers. We can actually lift $\h\phi_\Q$ globally, meaning we will not have to work section-by-section. Consider post-composing the functors $\O^\top$ and $\O$ with the localisation functor $\CAlg\to \CAlg_\Q$, and denote the resulting presheaves with a subscript $\Q$. By construction (also see \cite[Pr.4.47]{hilllawson}) the functor $\O^\top_\Q$ is formal, and by \cite[Pr.4.8]{lennartconnective} the sheaf $\O_\Q$ is also formal. This yields the following chain of equivalences lifting $\h\phi_\Q$:
\[\phi_\Q\colon \O^\top_\Q\xrightarrow{\simeq} \pi_\ast \O^\top_\Q\xrightarrow{\h\phi_\ast^\Q,\simeq} \pi_\ast\O_\Q\xleftarrow{\simeq} \O_\Q\]

%%%%%%%%%%%%%

\paragraph{(Transchromatic compatibility)}
We now have morphisms fitting into the following not a priori commutative solid diagram of presheaves of $p$-complete $\E_\infty$-rings on $\U$:
\begin{equation}\label{chromaticfracturecubetwo}\begin{tikzcd}
{\O^\top_p}\ar[dd]\ar[rr]\ar[dr, dashed, "{\phi_p}"]	&&	{\O^\top_{K(2)}}\ar[dd]\ar[rd, "{\phi_{K(2)}}"]	&	\\
	&	{\O_p}\ar[rr, crossing over]			&&	{\O_{K(2)}}\ar[dd]	\\
{\O^\top_{K(1)}}\ar[dr, "{\phi_{K(1)}}", swap]\ar[rr, "{\al^\top_\chrom\qquad\qquad}"]	&&	{(\O^\top_{K(2)})_{K(1)}}\ar[rd, "{(\phi_{K(2)})_{K(1)}}", swap]	&	\\
	&	{\O_{K(1)}}\ar[from=uu, crossing over]\ar[rr, "{\al_\chrom}", swap]				&&	{(\O_{K(2)})_{K(1)},}
\end{tikzcd}\end{equation}
The right face commutes from the naturality of the unit of the $K(1)$-localisation functor. We also claim that the bottom face commutes. In other words, we claim that for each $\Spec R\to \M$ inside $\U$, there is a natural path $\ga(R)$ between $\al_\chrom\circ \phi_{K(1)}$ and $(\phi_{K(2)})_{K(1)}\circ \al_\chrom^\top$ as maps of $\E_\infty$-$\tmf$-algebras. Note that the $\E_\infty$-$\tmf$-algebra structure on $(\O_{K(2)})_{K(1)}$ can come from either one of these maps (and a posteriori these two choices will agree up to homotopy). By \cite[p.44]{marktmf}, the $p$-adic $K$-theory functor induces an equivalence
\begin{equation}\label{discretepsaces}\Map_{\CAlg_{\tmf_{K(1)}}}(\EE^\top_{K(1)}, (\EE_{K(2)})_{K(1)})\xrightarrow{\simeq} \Hom_{\theta\Alg_{(V^\wedge_\infty)_\ast}}(\KU^\wedge_\ast \EE^\top_{K(1)}, \KU^\wedge_\ast (\EE_{K(2)})_{K(1)})).\end{equation}
As $\al_\chrom\circ\phi_{K(1)}$ and $(\phi_{K(1)})_{K(1)}\circ\al_\chrom^\top$ are isomorphic as functors into $\CAlg(\h\Sp)$ and hence also in the codomain of (\ref{discretepsaces}), we see these morphisms are homotopic as morphisms of $\E_\infty$-rings by the above equivalence. From the equivalence (\ref{discretepsaces}), we can employ \Cref{naturality} to obtain a homotopy between $\al_\chrom\circ\phi_{K(1)}$ and $(\phi_{K(2)})_{K(1)}\circ\al_\chrom^\top$ as morphisms of presheaves of $K(1)$-local $\E_\infty$-$\tmf_{K(1)}$-algebras from $\O^\top_{K(1)}$ to $(\O_{K(2)})_{K(1)}$. Using the fact that the front and back faces of (\ref{chromaticfracturecubetwo}) are Cartesian, we obtain a natural morphism of presheaves of $p$-complete $\E_\infty$-rings $\phi_p\colon \O_p^\top\to \O_p$ on $\U$, lifting $\h\phi_p$, as indicated by the dashed morphism in (\ref{chromaticfracturecubetwo}).

%%%%%%%%%%%%%

\paragraph{(Arithmetic compatibility)}
Currently, we have morphisms $\phi_\Q$ and $\phi_p$ fitting into the not a priori commutative solid diagram of presheaves of $\E_\infty$-rings on $\U$:
\begin{equation}\label{arithmeticfracturecubetwo}\begin{tikzcd}
{\O^\top}\ar[dd]\ar[rr]\ar[dr, dashed, "{\phi}"]	&&	{\prod_p \O_p^\top}\ar[dd]\ar[rd, "{\prod\phi_p}"]	&	\\
	&	{\O}\ar[rr, crossing over]			&&	{\prod_p \O_p}\ar[dd]	\\
{\O^\top_\Q}\ar[dr, "{\phi_\Q}", swap]\ar[rr]	&&	{\left(\prod_p \O_p^\top\right)_\Q}\ar[rd, "{(\prod \phi_p)_\Q}", swap]	&	\\
	&	{\O_\Q}\ar[from=uu, crossing over]\ar[rr]				&&	{\left(\prod_p \O_p\right)_\Q.}
\end{tikzcd}\end{equation}
Similar to the transchromatic compatibilities, the right face naturally commutes, so we are left to argue why the bottom face commutes. To study the bottom face, let us first work on the open substacks $\M_{\Ell}[c_4^{-1}]$ and $\M_{\Ell}[\Delta^{-1}]$ of $\M_{\Ell}$, which themselves form a cover of $\M_{\Ell}$; see \cite[\textsection9]{marktmf}. We then follow an analogous argument to the transchromatic situation above; see \cite[p.51]{marktmf} which shows the discreteness of the desired mappings spaces. Indeed, as the two homotopies witnessing the commutativity of the bottom face of (\ref{arithmeticfracturecubetwo}) restricted to the substacks $\M_{\Ell}[c_4^{-1}]$ and $\M_{\Ell}[\Delta^{-1}]$ agree on their intersection $\M_{\Ell}[c_4^{-1},\Delta^{-1}]$ (as the mapping spaces in question are discrete) these homotopies then glue to a homotopy on $\M_{\Ell}$. This yields a homotopy witnessing the commutativity of the bottom face of (\ref{arithmeticfracturecubetwo}). As the front and back faces of (\ref{arithmeticfracturecubetwo}) are Cartesian, we obtain our final natural equivalence of presheaves of $\E_\infty$-rings $\phi\colon \EE^\top\to \EE$ on $\U$, lifting $\h\phi$.\\

Therefore, $\Z$ is connected. The same argument can be made for $\Z^\sm$; see \Cref{simplification} for why \Cref{thetaalgebra} simplifies in this case. \qed

%%%%%%%%%%%%%%%%%%%%%%%%%%%%%%%%%%%%%%%%%%%%%%%%%%%%

\section{Compatibility of $\theta$-algebra structures}\label{finalcompsection}
The above proof of \Cref{mainguy} is contingent on \Cref{thetaalgebra}, whose proof we find rather delicate. Recall from \cite[\textsection6]{marktmf} that the $p$-adic $K$-theory of an $\E_\infty$-ring has the structure of a $\theta$-algebra, and this structure is functorial in morphisms of $\E_\infty$-rings.

\begin{lemma}\label{thetaalgebra}
Fix a prime $p$. Let $\O$ be an object of $\Z$ and $\h\phi\colon \h\O^\top\xrightarrow{\simeq}\h\O$ be the given equivalence of diagrams of homotopy commutative ring spectra. Then for any \'{e}tale $\Spec R\to \M_\Ell$, the map induced by
\[\h\phi\colon \FF^\top=\O_{K(1)}^\top(R)\to \O_{K(1)}(R)=\FF\]
on the zeroth $p$-adic $K$-theory ring is a morphism of $\theta$-algebras.
\end{lemma}

In general, it is not true that a morphism of homotopy commutative ring spectra should induce a morphism of $\theta$-algebras upon taking their $p$-adic $K$-theory, even if the homotopy commutative ring spectra involved come equipped with some $\E_\infty$-structures. \\

However, in the situation above the sections of the $K(2)$-localisation of the sheaf of $\E_\infty$-rings $\O$ have a prescribed $\E_\infty$-structure given by Lubin--Tate spectra (also called Morava $E$-theory); see the $K(2)$-local case in the proof of \Cref{mainguy} above. The comparison map in the chromatic fracture square between the $K(1)$-localisation of $\O$ and the $K(1)$-localisation of its $K(2)$-localisation is a map of $\E_\infty$-rings, and if we can show it induces an injection on $p$-adic $K$-theory we would obtain \Cref{thetaalgebra}. This is first done for an explicit \'{e}tale morphism into $\M_\Ell$, which has the properties that it covers $\M_\Ell^\sm$ and each of its connected component is an integral domain. Some descent and deformation theory to obtain this result for a general \'{e}tale morphism.

\begin{proof}
To show $\lambda\colon \K^\wedge_0 \FF^\top\to \K^\wedge_0 \FF$, the map induced by $\h\phi$ on $p$-adic $K$-theory, is a morphism of $\theta$-algebras, one must check it commutes with the stable $p$-adic Adams operations $\psi^\ell$ for every $\ell\in\mathbf{Z}_p^\times$ as well as the action of the operator $\theta$. The stable $p$-adic Adams operations $\psi^\ell$ are constructed on the spectrum $\KU_p$, so we automatically have compatibility with them for any map of spectra. It will be shown shortly that both rings above are \'{e}tale over the ring $V^\wedge_\infty$, hence they are $V^\wedge_\infty$-torsion free. In particular, this implies that both $\K^\wedge_0\FF^\top$ and $\K^\wedge_0\FF$ are $\mathbf{Z}_p$-torsion free, in which case the operator $\theta$ is equivalent datum to the $p$-adic Adams operator $\psi^p$; see \cite[Rmk.2.2.5]{gh04}. It therefore suffices show that the following diagram of $\mathbf{Z}_p$-algebras commutes:
\begin{equation}\label{firstdiagramhehe}\begin{tikzcd}
{\K^\wedge_0 \FF^\top}\ar[d, "{\psi^p_\top}"]\ar[r, "{\lambda}"]	&	{\K^\wedge_0 \FF}\ar[d, "{\psi^p_\top}"]	\\
{\K^\wedge_0 \FF^\top}\ar[r, "{\lambda}"]					&	{\K^\wedge_0 \FF}
\end{tikzcd}.\end{equation}
Let us write $R^\ord$ for the base-change of $\Spec R\to \M_\Ell$ over $\M^\ord_\Ell\to \widehat{\M}_\Ell\to \M_\Ell$, where $\M^\ord_\Ell$ is the moduli stack of generalised elliptic curves over $p$-complete rings whose reduction modulo $p$ is ordinary. By \cite[Pr.7.16]{marktmf}, we see $\FF^\top_{K(1)}$ is an elliptic cohomology theory for $\Spf R^\ord\to \M^\ord_\Ell$, and we can also consider $\FF_{K(1)}$ as an elliptic cohomology theory using $\h\phi_{K(1)}$. Define $W$ using the Cartesian diagram of formal stacks
\[\begin{tikzcd}
{\Spf W}\ar[d]\ar[r]	&	{\M^\ord_\Ell(p^\infty)}\ar[d]	\\
{\Spf R^\ord}\ar[r]		&	{\M^\ord_\Ell,}
\end{tikzcd}\]
where $\M^\ord_\Ell(p^\infty)$ is the formal stack of generalised elliptic curves $E$ over $\Spf R$ with ordinary reduction modulo $p$ with a given isomorphism $\eta\colon \widehat{\G}_m\to \widehat{E}$; see \cite[\textsection5]{marktmf}. The stack $\M^\ord_\Ell(p^\infty)$ is represented by the formal affine scheme $\Spf V_\infty^\wedge$ which is ind-\'{e}tale over $\M_\Ell^\ord$; see \cite[p.14-5]{marktmf}. This $W$ also has the structure of a $\theta$-algebra (see \cite[\textsection6]{marktmf}), and we denote the $p$-adic Adams operation on $W$ by $\psi^p_\alg$. By \cite[Pr.6.1]{marktmf}, or rather its proof, we obtain isomorphisms of $\mathbf{Z}_p$-algebras $v^\top\colon \K^\wedge_0 \FF^\top\simeq W$ and $v\colon \K^\wedge_0 \FF\simeq W$, which are natural in complex orientation preserving morphisms in $\CAlg(\h\Sp)$. These isomorphisms are \textbf{not} a priori isomorphisms of $\theta$-algebras; see \cite[\textsection6.2]{marktmf}. As $\FF$ obtains the structure of an elliptic cohomology theory for $R^\ord$ from the equivalence $\h\phi_{K(1)}$, we see that the following diagram of isomorphisms of $\mathbf{Z}_p$-algebras commutes:
\[\begin{tikzcd}
{\K^\wedge_0 \FF^\top}\ar[rr, "\lambda"]\ar[rd, swap, "{v^\top}"]	&&	{\K^\wedge_0 \FF}\ar[dl, "{v}"]	\\
	&	{W}	&
\end{tikzcd}\]
By construction (see \cite[Rmk.6.3]{marktmf}), we see $v^\top$ is an isomorphism of $\theta$-algebras. To show $\lambda$ is a morphism of $\theta$-algebras, it suffices to show $v$ is a morphism of $\theta$-algebras, or in other words: (\ref{firstdiagramhehe}) commutes if and only if the following diagram of $\mathbf{Z}_p$-algebras commutes:
\begin{equation}\label{secondlittleidagram}\begin{tikzcd}
{\K^\wedge_0 \FF}\ar[d, "{\psi^p_\top}"]\ar[r, "{v}"]	&	{W}\ar[d, "{\psi^p_\alg}"]	\\
{\K^\wedge_0 \FF}\ar[r, "{v}"]					&	{W}
\end{tikzcd}\end{equation}
Let us now prove this is the case for a specific \'{e}tale map $\Spec R\to \M_\Ell$.

\paragraph{(Choosing a particular \'{e}tale morphism)}
Recall the moduli stack $\M^\sm_1(N)$ of smooth elliptic curves with $\Ga_1(N)$-level structure, from \cite[(1.3.12)]{handbooktmf} for example. Importantly, recall the map $\M^\sm_1(N)\to \M_{\Ell, \mathbf{Z}[\frac{1}{N}]}^\sm$ is an \'{e}tale cover and that for $N\geq 4$ the moduli stack $\M^\sm_1(N)$ is in fact affine. This implies that the morphism of stacks
\[\Spec A=\M_1(4)\sqcup \M_1(5)\to \M_{\Ell,\mathbf{Z}[\frac{1}{2}]}^\sm\sqcup\M_{\Ell,\mathbf{Z}[\frac{1}{5}]}^\sm\to \M_\Ell^\sm\to \M_\Ell\]
is \'{e}tale, and the restriction to $\M_\Ell^\sm$ is an \'{e}tale cover. By base-change over $\Spf \mathbf{Z}_p$ we obtain an \'{e}tale map $E\colon \Spf \widehat{A}\to \widehat{\M}_\Ell$. Following \cite[p.42-3]{marktmf}, write $A_\ast$ for the graded ring defined by $A_{2k}=\omega^{\otimes k}_E(\Spf \widehat{A})$ were $E$ is the elliptic curve over $\widehat{A}$ defined by the map of formal stacks above. Note that the Hasse invariant $v_1$ for $E$ lives in $A_{2(p-1)}$. Let us also make the following definitions:
\begin{equation}\label{variousrings}A^\ord_\ast=(A_\ast)[v_1^{-1}]^\wedge_p,\qquad A^\supers_\ast=(A_\ast)^\wedge_{(v_1)},\qquad (A_\ast^\supers)^\ord=(A_\ast^\supers)[v_1^{-1}]^\wedge_p.\end{equation}
If we omit the subscript $\ast$ we are implicitly considering the ring in degree zero. By \cite[(8.6)]{marktmf}, there is a canonical map $\al_\ast\colon A_\ast^\ord\to (A_\ast^\supers)^\ord$ as $v_1$ is invertible in $(A_\ast^\supers)^\ord$, and we now define $W^\supers$ using the diagram of stacks
\begin{equation}\label{hereisdiagramheremeroar}\begin{tikzcd}
{\Spf W^\supers}\ar[r, "{\widetilde{\al}}"]\ar[d]	&	{\Spf W}\ar[d]\ar[r]	&	{\M^\ord_\Ell(p^\infty)}\ar[d]	\\
{\Spf (A^\supers)^\ord}\ar[r, "{\al}"]		&	{\Spf A^\ord}\ar[r]		&	{\M^\ord_\Ell,}
\end{tikzcd}\end{equation}
where all squares are Cartesian. The ring $W^\supers$ obtains a $\theta$-algebra structure from the above diagram, and in such a way that $\widetilde{\al}\colon W\to W^\supers$ is a morphism of $\theta$-algebras; see \cite[p.40]{marktmf}. We claim that $\widetilde{\al}$ comes from a map of $\E_\infty$-rings.

\begin{claim}\label{claimone}
The zeroth $p$-adic $K$-theory of the canonical map of $\E_\infty$-rings
\[\al_\chrom\colon \FF^\ord=\O_{K(1)}(A)\to (\O_{K(2)}(A))_{K(1)}=(\FF^\supers)^\ord\]
is isomorphic to $\widetilde{\al}$.
\end{claim}

\begin{proofclaim}{Proof of \Cref{claimone}}
We have already seen that $\FF^\ord=\FF_{K(1)}$ is an elliptic cohomology theory for the map $\Spf A^\ord\to \M^\ord_\Ell$, and similarly by \cite[Lm.8.8]{marktmf}, we see that $(\FF^\supers)^\ord$ is an elliptic cohomology theory for the map $\Spf (A^\supers)^\ord\to \M^\ord_\Ell$. The same is true for $\A=\O^\top_{K(1)}(A)$ and $(\O^\top_{K(2)}(A))_{K(1)}=\A'$, and in this case we know that taking $\pi_0$ of $\al_\chrom^\top\colon \A\to \A'$ is isomorphic to $\al\colon A^\ord\to (A^\supers)^\ord$ by construction; see \cite[p.43-4]{marktmf}. The naturality of $\h\phi\colon \O^\top\to\O$ and the chromatic fracture square imply that $\pi_0$ of the natural map of $\E_\infty$-rings $\al_\chrom\colon \FF^\ord\to (\FF^\supers)^\ord$ also realises $\al$, and hence taking zeroth $p$-adic $K$-theory realises $\widetilde{\al}$. This proves \Cref{claimone}.
\end{proofclaim}

Recall that $\FF^\ord=\O_{K(1)}(A)$ for our choice of $A$ above. Consider the diagram of $\mathbf{Z}_p$-algebras
\begin{equation}\label{coop}\begin{tikzcd}
{\K^\wedge_0 \FF^\ord}\ar[dd]\ar[rr]\ar[dr]	&&	{W}\ar[dd]\ar[rd]	&	\\
	&	{\K^\wedge_0 (\FF^\supers)^\ord}\ar[rr, crossing over]			&&	{W^\supers}\ar[dd]	\\
{\K^\wedge_0\FF^\ord}\ar[dr]\ar[rr]		&&	{W}\ar[rd, "{\widetilde{\al}}"]		&	\\
	&	{\K^\wedge_0 (\FF^\supers)^\ord}\ar[from=uu, crossing over]\ar[rr]	&&	{W^\supers,}
\end{tikzcd}\end{equation}
where the maps are the obvious ones used above, and all the vertical morphisms are the unstable Adams operations; $\psi^p_\top$ on the left, and $\psi^p_\alg$ on the right. Note that the top and bottom faces commute by \Cref{claimone}, the right face commutes as $\widetilde{\al}\colon W\to W^\supers$ is a morphism of $\theta$-algebras, and the left face commutes $\al_\chrom$ is a morphism of $\E_\infty$-rings. Most importantly, the front face also commutes. Indeed, from the arguments in the $K(2)$-local case of the proof of \Cref{mainguy} we see $\FF^\supers$ is naturally equivalent to a product of $K(2)$-local Lubin--Tate spectra recognising the given elliptic curve over $\Spf A^\supers$, and we can then apply \cite[Th.6.10]{marktmf}; the hypotheses and proof of this theorem are dispersed between pages 21 and 24 of \emph{idem}. The back face of (\ref{coop}) is precisely (\ref{secondlittleidagram}) for $R=A$.

\begin{claim}\label{injectivity}
The morphism $\widetilde{\al}\colon W\to W^\supers$ is injective.
\end{claim}

Using this claim for now, to show that the back face of (\ref{coop}) commutes, it suffices to do so after post-composing with $\widetilde{\al}$. This follows from the above considerations by a diagram chase. Hence the back face of (\ref{coop}) commutes, which yields the commutativity of (\ref{secondlittleidagram}) for this particular choice of \'{e}tale map $\Spec A\to \M_\Ell$.

\begin{proofclaim}{Proof of \Cref{injectivity}}
As $\M^\ord_\Ell(p^\infty)\to \M^\ord_\Ell$ is ind-\'{e}tale (see \cite[Lm.5.1]{marktmf}), it is flat. By base-change, we see that $A^\ord\to W$ is also flat, hence $\widetilde{\al}$ is injective if we can show $\al$ is injective. To do this, we will show $\al_\ast\colon A_\ast^\ord\to (A_\ast^\supers)^\ord$ is injective. Using the notation above, we find ourselves with the following commutative diagram of graded rings, where all maps are the indicated localisations or completions:
\[\begin{tikzcd}
{A_\ast}\ar[r]\ar[d, "{\ga}"]	&	{A_\ast[v_1^{-1}]}\ar[d, "{\be}"]\ar[r]	&	{A_\ast[v_1^{-1}]^\wedge_p=A_\ast^\ord}\ar[d, "{\al_\ast}"]	\\
{A_\ast^\supers}\ar[r]			&	{A_\ast^\supers[v_1^{-1}]}\ar[r]		&	{A_\ast^\supers[v_1^{-1}]^\wedge_p=(A^\supers_\ast)^\ord.}
\end{tikzcd}\]
Let us now make the following remarks from this diagram:
\begin{enumerate}
\item From our choice of $A$, we have $A=A_1\times A_2$, where $A_1$ and $A_2$ both integral domains; see \cite[Th.1.1.1]{markkyle} for the $\M^\sm_1(5)$-case, and the $\M^\sm_1(4)$-case is similar.\footnote{Indeed, following the proof of \cite[Th.1.1.1]{markkyle}, which in turn uses \cite[\textsection4.4]{husemollerell}, each elliptic curve $E$ in Weierstra\ss form over a $\mathbf{Z}[\frac{1}{2}]$-algebra $R$ with $(0,0)$ a point of order $4$ can be moved into (nonhomogeneous) \emph{Tate normal form} $y^2+(1-c)xy-by=x^3-bx^2$. As $(0,0)$ has order $4$ we have $[2](0,0)=[-2](0,0)$, which by \cite[Ex.4.4]{husemollerell} yields $c=0$. As in the proof of \cite[Th.1.1]{markkyle}, we see $\M_1(4)$ is equivalent to spectrum of $\mathbf{Z}[\frac{1}{2}, b, \Delta^{-1}]$ where $\Delta=b^4(1+16b)$.} It then follows that $\ga$ can be written in the following commutative diagram of graded rings
\[\begin{tikzcd}
{A_\ast}\ar[r, "{\ga}"]\ar[d]						&	{A_\ast^\supers=(A_\ast)^\wedge_{(v_1)}}\ar[d]	\\
{{A_{\ast,1}}\times {A_{\ast,2}}}\ar[r, "{\ga_1\times\ga_2}"]	&	{A_{\ast,1}^\supers\times A_{\ast,2}^\supers.}
\end{tikzcd}\]
The ring $A$ is Noetherian as it is finitely presented over $\Spec \mathbf{Z}$, so both $A_1$ and $A_2$ are Noetherian integral domains. In particular, the completion maps $\ga_i$ are flat for $i=1,2$. If we know these maps $\ga_i$ are nonzero, then it immediately follows that they are injective. To see that they are nonzero, it suffices to show that $v_1$ is not a unit inside both $A_{\ast,1}$ and $A_{\ast,2}$. This is where our choice of $A$ comes in. If our fixed prime $p\neq 2,5$, then for both $i=1,2$ the image of the map
\[\Spf \widehat{A}_i\to \widehat{\M}_\Ell^\sm\to \widehat{\M}_\Ell\]
contains a supersingular elliptic curve, as all supersingular elliptic curves are contained in the smooth locus of $\widehat{\M}_\Ell$. This implies that $v_1$ cannot be a unit, else $\Spf \widehat{A}_i\to \widehat{\M}_\Ell$ would define only ordinary elliptic curves of height one. Similarly, if $p=2$, then the $p$-completion of $A$ is $\widehat{A}_2$, and we again see $v_1$ is not a unit so $\ga_2=\ga$ is injective. The same holds for $\widehat{A}_1$ when $p=5$. This implies that $\ga_1\times\ga_2$ is always injective, hence $\ga$ is injective.
\item As $\be$ is the $v_1$-localisation of $\ga$, and localisation is exact, we see that $\be$ is also injective.
\item Standard arguments show that the $p$-completion of $\be$, also known as $\al_\ast$, is also injective. Indeed, limits are left exact, so it suffices to show each $\al_\ast^k$ in the following commutative diagram of rings is injective, for every $k\geq 1$:
\[\begin{tikzcd}
{A_\ast[v_1^{-1}]}\ar[r]\ar[d, "{\be}"]	&	{A_\ast[v_1^{-1}]/p^k}\ar[d, "{\al_\ast^k}"]	\\
{A^\supers_\ast[v_1^{-1}]}\ar[r]		&	{A^\supers_\ast[v_1^{-1}]/p^k.}
\end{tikzcd}\]
Given an element $\overline{x}$ such that $\al_\ast^k(\overline{x})=0$, then we first note that any lift $x$ over $\overline{x}$ is sent to a $\be(x)$ such that $p^k\be(x)=0$. However, $A^\supers_\ast[v_1^{-1}]$ is flat over $\mathbf{Z}$, as we have the following composite of flat maps:
\[\mathbf{Z}\to \mathbf{Z}_p\to \widehat{A} \to A_\ast \xrightarrow{\ga} A_\ast^\supers\to A^\supers_\ast[v_1^{-1}];\]
the second map is flat as $\Spec A\to \M_\Ell$ is \'{e}tale and $\M_\Ell$ is smooth over $\mathbf{Z}$, and the third map is flat as each $\omega^{\otimes k}_E(\Spf A)$ is a line bundle and hence projective of rank 1. This implies that $A^\supers_\ast[v_1^{-1}]$ is torsion-free, hence $\be(x)=0$. As $\be$ is injective, this implies $x=0$ and $\overline{x}=0$, hence $\al_\ast^k$ is injective.
\end{enumerate}
It follows that $\widetilde{\al}$ is injective.
\end{proofclaim}

\paragraph{Reduction to a general \'{e}tale morphism}
Let $\Spec R\to \M_\Ell$ be an arbitray \'{e}tale morphism now, and consider the Cartesian diagram of stacks
\[\begin{tikzcd}
{\Spec B}\ar[r]\ar[d]	&	{\Spec A}\ar[d]	\\
{\Spec R}\ar[r]		&	{\M_\Ell,}
\end{tikzcd}\]
where $\Spec A=\M_1^\sm(4)\sqcup\M_1^\sm(5)$ is that of the previous paragraph; the stack $\Spec B$ is affine as $\M_\Ell$ is separated. All of the morphisms above are \'{e}tale by base-change, so we can consider the morphism of $\E_\infty$-rings $\O(A)\to \O(B)$.

\begin{claim}\label{etaleclaim}
The morphism of $\E_\infty$-rings $\O(A)\to \O(B)$ is \'{e}tale.
\end{claim}

\begin{proofclaim}{Proof}
Recall from \cite[\textsection7.5]{haname} that a morphism $\A\to \B$ of $\E_\infty$-rings is \emph{\'{e}tale} if the morphism $\pi_0 \A\to \pi_0 \B$ of discrete commutative rings is \'{e}tale and the natural map of $\pi_0 \B$-modules
\[\pi_0\B\underset{\pi_0\A}{\otimes} \pi_\ast \A\to \pi_\ast\B\]
is an isomorphism. The fact that $\pi_0 \O(A)\to \pi_0 \O(B)$ is \'{e}tale follows from the facts that $A\to B$ is \'{e}tale and $\O$ defines a natural elliptic cohomology theory. The condition on the higher homotopy groups also follows as $\O$ defines a natural elliptic cohomology theory; see \Cref{higherhomotopy}.
\end{proofclaim}

By \cite[Th.7.5.0.6]{haname}, the $\pi_0$-functor induces an equivalence of $\infty$-categories
\[\CAlg^\et_{\O(A)}\xrightarrow{\pi_0} \CAlg^\et_A,\]
where the superscript indicates subcategories of \'{e}tale algebras. By \Cref{etaleclaim}, for any \'{e}tale $\E_\infty$-$\O(A)$-algebra $\B$ such that $\pi_0 \B$ is isomorphic to $B$ as an $A$-algebra, there is a equivalence of $\E_\infty$-$\O(A)$-algebras $\O(B)\simeq \B$, which is unique up to contractible choice. As we have proven \Cref{thetaalgebra} for $\Spec A$, it follows from the proof of \Cref{mainguy} above that the equivalence of homotopy commutative ring spectra $\h\phi(A)\colon \O^\top(A)\simeq \O(A)$ can be lifted to a morphism of $\E_\infty$-rings. The composite $\O(A)\simeq \O^\top(A)\to \O^\top(B)$ is also an \'{e}tale $\E_\infty$-$\O(A)$-algebra recognising $B$, hence we obtain a natural equivalence of $\E_\infty$-$\O(A)$-algebras $\O^\top(B)\simeq \O(B)$. As $\O^\top(B)$ is $\theta$-compatible (see \cite[Rmk.6.3]{marktmf}), we see $\O(B)$ is also $\theta$-compatible, meaning that (\ref{secondlittleidagram}) commutes for $R=B$. Finally, let us turn our attention to $\O(R)\to \O(B)$.

\begin{claim}\label{finalinjectivitythingy}
The morphism induced by $\O(R)\to \O(B)$ on zeroth $p$-adic $K$-theory is injective.
\end{claim}

Assuming the above claim, it immediately follows that (\ref{secondlittleidagram}) for our arbitrary $R$. Indeed, \Cref{finalinjectivitythingy} provides us with an injection of $\theta$-algebras induced by $\O(R)\to \O(B)$, which allows us to check the commutativity of (\ref{secondlittleidagram}) in the same diagram for $R=B$, which we know commutes by the above paragraph.

\begin{proofclaim}{Proof of \Cref{finalinjectivitythingy}}
The morphism $\Spec B\to \Spec R$ can be factored into the following diagram of formal stacks:
\begin{equation}\label{bigdiagram}\begin{tikzcd}
{\Spf W_B}\ar[r]\ar[d]			&	{\Spf W_R^\sm}\ar[d]\ar[r]			&	{\Spf W_R}\ar[r]\ar[d]			&	{\M^\ord_\Ell(p^\infty)}\ar[d]	\\
{\Spf \widehat{B}^\ord}\ar[r]\ar[d]	&	{\Spf \widehat{R}^{\ord,\sm}}\ar[r]\ar[d]	&	{\widehat{R}^\ord}\ar[r]\ar[d]		&	{\M^\ord_\Ell}\ar[d]		\\
{\Spf \widehat{B}}\ar[r]\ar[d]		&	{\Spf \widehat{R}^\sm}\ar[r]\ar[d]		&	{\Spf\widehat{R}}\ar[r]\ar[d]		&	{\widehat{\M}_\Ell}\ar[d]		\\
{\Spec B}\ar[d]\ar[r]			&	{\Spec R^\sm}\ar[d]\ar[r]				&	{\Spec R}\ar[r]\ar[d]			&	{\M_\Ell}				\\
{\Spec A}\ar[r]				&	{\M_\Ell^\sm}\ar[r]					&	{\M_\Ell}
\end{tikzcd}\end{equation}
Every square above is Cartesian, and the $\widehat{(-)}$ indicates base-change with $\Spf\mathbf{Z}_p$. By \cite[Lm.6.1]{marktmf}, the morphism $W_R\to W_B$ above is isomorphic to the morphism induced by $\O(R)\to\O(B)$ on $p$-adic $K$-theory, hence it suffices to see the composite map
\begin{equation}\label{littlecomposition}W_R\to W_R^\sm\to W_B,\end{equation}
featured in the top-left corner of (\ref{bigdiagram}), is injective. As $\Spec A\to \M_\Ell^\sm$ is an \'{e}tale cover, then by base-change we see $W^\sm_R\to W_B$ is also faithfully flat, and hence injective. Observe that $W_R\to W_R^\sm$ is an open immersion of formal affine schemes by base-change as $\M^\sm_\Ell\to \M_\Ell$ is an open immersion of stacks. Moreover, we claim the open immersion $R\to R^\sm$ has scheme theoretically dense image as $\Delta$ is a nonzero divisor in $R$; see \cite[\href{https://stacks.math.columbia.edu/tag/01RE}{Tag 01RE}]{stacks}. Indeed, to see $\Delta$ is not a zero divisor, it suffices to show that the image of $\Spec R\to \M_\Ell$ has nontrivial intersection with the image of $\M_\Ell^\sm$. This is clear on the level of underlying topological spaces, as the inclusion $|\M_\Ell^\sm|\to |\M_\Ell|$ is equivalent to open immersion of coarse moduli spaces $|\mathbf{A}^1_\mathbf{Z}|\to |\mathbf{P}^1_\mathbf{Z}|$ which adds the point at $\infty$, and the map $|\Spec R|\to |\M_\Ell|$ is open as \'{e}tale morphisms are in particular flat and locally of finite presentation; see \cite[\href{https://stacks.math.columbia.edu/tag/06R7}{Tag 06R7}]{stacks}. As all the right vertical maps in (\ref{bigdiagram}) are flat, and $R\to R^\sm$ is quasi-compact (as a map of affine schemes), then \cite[\href{https://stacks.math.columbia.edu/tag/0CMK}{Tag 0CMK}]{stacks} tells us that $W_R\to W_R^\sm$ also have scheme theoretically dense image. Another application of \cite[\href{https://stacks.math.columbia.edu/tag/01RE}{Tag 01RE}]{stacks} show this open immersion $W_R\to W_R^\sm$ must be injective. Therefore, the composite (\ref{littlecomposition}) is injective.
\end{proofclaim}

This finishes our proof of \Cref{thetaalgebra}.
\end{proof}

\begin{remark}\label{simplification}
Let us note a potential improvement that could be made to this note. If the reader can find an affine \'{e}tale \textbf{cover} $\Spec A\to \M_\Ell$ such that \Cref{injectivity} holds, which is a purely algebro-geometric pursuit, then the rest of \Cref{mainguy} follows. Indeed, in this case one can prove that $\O^\top$ is uniquely defined on the \v{C}ech nerve of such a cover by copying the proof of \Cref{mainguy} and \Cref{thetaalgebra} seen above. One can then use spectral deformation theory and descent to show that \Cref{mainguy} follows from this particular case. With this in mind, the reader might also notice that restricting \Cref{tha} to the moduli stack of smooth elliptic curves vastly simplified from the above proof.
\end{remark}

%REFERENCES!!!
\addcontentsline{toc}{section}{References}
%bib on laptop %%uncomment below at home!%%
\bibliography{C:/Users/jackd/Dropbox/Work/references} 

%bib at work %%uncomment below at work!%%
%\bibliography{/Users/davie006/Desktop/Dropbox/Work/references} 

\bibliographystyle{alpha}

%\newpage
%\listoftodos

\end{document}